\documentclass[12pt]{article}
\usepackage[T2A]{fontenc}
\usepackage[cp1251]{inputenc}
\usepackage[russian]{babel}
\usepackage{amssymb,amsmath,latexsym,amsthm}
\usepackage{amsbsy,amsfonts,amssymb,mathrsfs,amscd}
\usepackage{eucal}
\usepackage{graphicx}
\usepackage{indentfirst}
\usepackage{verbatim}
\usepackage{enumerate}
\usepackage[pdftex,unicode,colorlinks,linkcolor=blue,
citecolor=red,bookmarksopen,pdfhighlight=/N]{hyperref}

\theoremstyle{plain}

\newtheorem{theorem}{Теорема}

\theoremstyle{definition}
\newtheorem*{definition}{Определения}
\newtheorem*{example}{Пример}

\newcommand{\CC}{\mathbb C}
\newcommand{\NN}{\mathbb N}
\newcommand{\RR}{\mathbb R}
\newcommand{\BB}{\mathbb B}
\newcommand{\dd}{\,{\rm d}}
\newcommand{\const}{{\rm const}}

\renewcommand{\leq}{\leqslant}
\renewcommand{\geq}{\geqslant}

\DeclareMathOperator{\Hol}{Hol}
\DeclareMathOperator{\sbh}{sbh}
\DeclareMathOperator{\Int}{int}
\DeclareMathOperator{\clos}{clos}
\DeclareMathOperator{\Zero}{Zero}
\DeclareMathOperator{\supp}{supp}

\title{\begin{flushleft} {\normalsize УДК 517.55 + 517.574 + 517.987.1 + 517.53}\end{flushleft} {\large \bf К распределению нулевых множеств голоморфных функций. II\footnote{Исследование выполнено за счёт гранта Российского научного фонда (проект № 18-11-00002).}}}
\author{Б.\,Н.~Хабибуллин, Э.\,Б.~Хабибуллина}
\date{10 июня 2018 г.}
\begin{document}
\maketitle 
Обобщается [1, теорема 3] в обозначениях и определениях из [1]. Для $S$ и $\mathcal O$ из $\CC^n_{\infty} :=\CC^n\cup\{\infty\}$, $n\in \NN$,  {\it замыкание\/} $\clos S$, {\it внутренность\/} $\Int S$ и {\it граница\/} $\partial S$ --- в $\CC^n_{\infty}$; $S\Subset \mathcal O$, если $\clos S \subset \mathcal O$; $\sbh(S)$ и $\Hol(S)$ ---  классы соответственно {\it субгармонических\/} и {\it голоморфных функций\/} в открытой окрестности $S$; $\sbh^+(S):= \{v\in \sbh (S)\colon v\geq 0\text{ на $S$}\}$. Всюду далее 
\begin{subequations}\label{DS}
\begin{align}
D&\text{ \it --- область в $\mathbb C^n_{\infty}\neq D$
и существует точка   $z_0 \in \Int S\Subset D$;}
\tag{\ref{DS}d}\label{DSD}
\\  
\sbh_0&(D\setminus S):=\Bigl\{v\in \sbh(D\setminus S)\colon \lim_{D\ni z'\to z}v(z')=0 \text{ \it   для всех\/  $z\in \partial D$}\Bigr\}.
\tag{\ref{DS}s}\label{DSs}
\end{align}
\end{subequations}

\begin{definition}[{\rm ср. с [1, определение 1]}]\label{def:tesf} 
Функцию $v\overset{\eqref{DSs}}{\in} \sbh_0 (D\setminus S)$ называем {\it вещественной  тестовой  функцией для\/ $D$ вне\/} $S$, если эта функция $v$ одновременно 
\begin{enumerate}[{\rm (a)}]
\item\label{wt:i}  {\it ограничена на\/ \/} $S_0\setminus S$ для некоторого  $S_0\subset D$, удовлетворяющего условию\/ $S\overset{\eqref{DSD}}{\Subset} \Int S_0$;   
 
\item\label{wt:ii}  {\it положительна,\/} т.\,е. $\geq 0$, 
{\it на\/} $D\setminus S_v$ для некоторого подмножества $S_v\Subset D$. 
\end{enumerate}
Класс всех таких функций $v$ обозначаем через  $\sbh_0^{\pm}(D\setminus S; <+\infty)$. В  [1, (2.2b)] определялся подкласс $\sbh_0^{+}(D\setminus S; <+\infty) :=\sbh_0^{\pm}(D\setminus S; <+\infty)\cap\, \sbh^+(D\setminus S)$ {\it положительных тестовых функций для\/ $D$ вне\/} $S$.  Для $b\in \RR^+:=\{x\in \RR\colon x\geq 0\}$ и $S\Subset \Int S_0\subset S_0\subset D$ полагаем
\begin{equation}\label{<b}
\sbh_0^{\pm}(D\setminus S,S_0;\leq b):=\Bigl\{v\in \sbh_0^{\pm} (D\setminus S;<+\infty)\colon \sup_{S_0\setminus S}|v|\leq b\Bigr\}\subset \sbh_0^{\pm}(D\setminus S; <+\infty).
\end{equation}
Подкласс 
$\sbh_0^{+}(D\setminus S;\leq b):=\sbh_0^{\pm}(D\setminus S,D;\leq b)\cap \,\sbh^+(D\setminus S)$ играл ключевую роль в [1].
Через $\sigma_{2n-2}$ обозначаем $(2n-2)$-{\it меру Хаусдорфа,\/} нормированную как в [1,1.2.3, (1.3)]. Для функции $f\in \Hol(D)$ через $\Zero_f$ обозначаем {\it дивизор нулей\/} функции $f$ [1,1.2.4]. 
\end{definition}

\begin{theorem}[{\rm обобщение [1, теорема 3]}]\label{th1} 
Пусть     $M(z_0)\neq -\infty$ для   $M\in \sbh (D)$ с мерой  Рисса $\nu_M$, $b\in \RR^+\setminus \{0\}$ и  $S\overset{\eqref{wt:i}}{\Subset} \Int S_0\subset S_0\subset D$. Тогда найдутся постоянные\footnote{Постоянная $\const^+_{a_1,a_2, \dots}\geq 0$ зависит только от $a_1,a_2,\dots$, а также от $n$ и $D$.} $C:=\const_{z_0,S,S_0,b}^+$
и $\overline C_M:=\const_{z_0,S, S_0, M}^+$, с которыми для каждой ненулевой функции $f\in{\Hol} (D)$ при ограничении $|f|\leq \exp M$ на $D$  для всех  $v\overset{\eqref{<b}}{\in}  \sbh_0^{\pm}(D\setminus S, S_0;\leq b)$ выполнено неравенство
\begin{equation}\label{iq:ufu}
\int_{D\setminus {S}}v\, {\Zero}_f  \dd \sigma_{2n-2}
{\leq} \int_{D\setminus {S}}  v \dd {\nu}_M	-C \log \bigl|f(z_0)\bigr| +C\,\overline C_M. 
\end{equation}
\end{theorem}
В [1, теорема 3] равномерное по  $v\in \sbh_0^{\pm}(D\setminus S, S_0;\leq b)$ неравенство \eqref{iq:ufu} дано  для {\it более узкого\/} (см. Пример ниже)  класса $\sbh_0^{+}(D\setminus S;\leq b)$ {\it положительных\/} тестовых функций. 

\begin{proof} Можно считать, что $z_0=0$ и  $u(0)=0$
для $u:=\log|f|\in \sbh (D)$ с мерой Рисса $\nu_u$. По условию 
\eqref{wt:ii} и определениям \eqref{DSs} и \eqref{<b}  функция $v \in \sbh_0^{\pm}(D\setminus S, S_0;\leq b)$, {\it продолженная нулем на\/ $\CC_{\infty}^n\setminus D$,  субгармонична на\/} $\CC^n_{\infty}\setminus S$. 
Существует регулярная для задачи Дирихле область $U\subset \CC^n_{\infty}$, удовлетворяющая условиям $S\Subset U\Subset \Int S_0$, с функцией Грина $g_U:=g_U(\cdot, 0) \geq 0$ с полюсом в нуле. После продолжения нулем на $\CC_{\infty}^n\setminus U$ функция $g_U$ субгармонична на $\CC^n_{\infty}\setminus \{0\}$ [2, 5.7.4]. Существует открытое множество $\mathcal O\Subset (\Int S_0)\setminus \clos S$, содержащее границу $\partial U\Subset \mathcal O$. Положим, следуя схеме [3, доказательство леммы 6.2],  
\begin{equation}\label{vv}
c:=\dfrac{2b+2}{\inf\limits_{U\cap \partial \mathcal O}g_U}>0, \qquad 
\widetilde{v}:=\begin{cases}
v \quad &\text{на $\CC_{\infty}^n\setminus U$},\\ 
\max\{v,cg_U-b-1\}\quad &\text{на $U\cap \mathcal O$},\\
cg_U-b-1 \quad &\text{на $U\setminus  \mathcal O$}.
\end{cases}
\end{equation}
Построение \eqref{vv} <<склеивает>> функцию $\widetilde{v}$ в субгармоническую на $\CC_{\infty}^n\setminus \{0\}$, тождественно равную нулю на $\CC_{\infty}^n\setminus D$ и  удовлетворяющую предельному соотношению 
\begin{equation}\label{vh}
\lim_{0\neq z\to 0} \frac{\widetilde{v}(z)}{-h_{2n}(|z|)}\overset{\eqref{vv}}{=} c, \quad\text{где }h_{2n}(t):=\begin{cases}
\log t\quad&\text{при $n=1$},\\
 -t^{2-2n}\quad&\text{при $n>1$},
\end{cases}
\quad t>0.
\end{equation} 
Рассмотрим теперь  {\it положительную вне некоторого\/} $S_v\overset{\eqref{wt:ii}}{\Subset} D$ функцию 
\begin{equation}\label{V}
\widetilde{V}\overset{\eqref{vv}}{:=}\frac{1}{c}\,\widetilde{v}\in \sbh(\CC_{\infty}^n\setminus \{0\}), \quad \lim_{0\neq z\to 0} \frac{\widetilde{V}(z)}{-h_{2n}(|z|)} \overset{\eqref{vh}}{=} 1, \quad \widetilde{V}\equiv 0 \text{ на $\CC^n_{\infty}\setminus D$}. 
\end{equation}
Каждому числу $\varepsilon >0$ сопоставим 
открытое множество $\mathcal O_{\varepsilon}=\{z\in \CC_{\infty}^n\setminus \{0\}\colon \widetilde{V}(z)<\varepsilon\}$ и функцию 
$V_{\varepsilon}$, тождественно равную нулю в каждой связной компоненте множества $\mathcal O_{\varepsilon}$, пересекающейся с $\CC_{\infty}^n\setminus D$, и равную функции $\widetilde{V}-\varepsilon$ в остальной части $\CC_{\infty}^n\setminus \{0\}$. Тогда при всех достаточно малых $\varepsilon>0$ в силу положительности $\widetilde{V}$ вне некоторого  $S_v\Subset D$ и ввиду условия $\lim_{D\ni z\to \partial  D} \widetilde{V}(z)\overset{\eqref{DSs}}{=}0$ 
функции $V_{\varepsilon}\in \sbh(\CC_{\infty}^n\setminus\{0\})$ согласно \eqref{V} тождественно равны нулю вне некоторого $S_{v,\varepsilon}\Subset D$ и удовлетворяют условию 
$\lim\limits_{0\neq z\to 0} \dfrac{V_{\varepsilon}(z)}{-h_{2n}(|z|)}
\overset{\eqref{V}}{=} 1$. Функции $V_{\varepsilon}$ с мерой Рисса $\mu_{\varepsilon}$, $\supp \mu_{\varepsilon}\Subset D\setminus \{0\}$, называем  
 {\it функциями,\/} или {\it потенциалами, Аренса\,--\,Зингера\/}, а $\mu_{\varepsilon}$ --- это {\it меры Аренса\,--\,Зингера\/} [4, определения 1, 2]. По [4, предложение 1.2] 
\begin{equation*}\label{uM}
\int_{D\setminus \{0\}} V_{\varepsilon} \dd {\nu}_u =\int_{D} u \dd  \mu_{\varepsilon}, 
\quad \int_{D} M \dd  \mu_{\varepsilon}=\int_{D\setminus \{0\}} V_{\varepsilon} \dd {\nu}_M+M(0).
\end{equation*}   
Отсюда, поскольку $u=\log |f|\leq M$, а по построению 
$V_{\varepsilon}\leq \widetilde{V}$ всюду на $\CC_{\infty}^n\setminus \{0\}$, получаем 
\begin{equation}\label{uMc}
\int_{D\setminus \{0\}} V_{\varepsilon} \dd {\nu}_u \leq \int_{D\setminus \{0\}} V_{\varepsilon} \dd {\nu}_M +M(0)
\leq \int_{D\setminus \{0\}} \widetilde{V} \dd {\nu}_M +M(0)
\quad\text{для всех $V_{\varepsilon}$}.
\end{equation}
По построению при $0<\varepsilon \to 0$ функции $V_{\varepsilon}\leq \widetilde{V}$ монотонно поточечно стремятся к $\widetilde{V}$. При $u=\log |f|$ по теореме о монотонной сходимости интегралов из \eqref{uMc} имеем
\begin{equation}\label{utof}
\int_{D\setminus \{0\}}\widetilde{V}\, {\Zero}_f  \dd \sigma_{2n-2} =\int_{D\setminus \{0\}} \widetilde{V} \dd {\nu}_u \overset{\eqref{uMc}}{\leq} \int_{D\setminus \{0\}} \widetilde{V} \dd {\nu}_M +M(0),
\end{equation}
где равенство в левой части следует из формулы Пуанкаре\,--\,Лелона [1, 1.2.4]. Для перехода от $\widetilde{V}$ к подынтегральной функции $\widetilde{v}$ необходимо домножить соотношение \eqref{utof} на число $c>0$ из \eqref{vv}--\eqref{V}. Затем, исходя из построения  $\widetilde{v}$ в \eqref{vv}, подобно [1, 3.3, после (3.26)], можно  избавиться от частей интегралов по множеству $S\Subset D$ и оставить в интегралах из  
\eqref{utof} только подынтегральную функцию $v$ вместо $\widetilde{v}$ из \eqref{vv}. Аккуратный контроль за возникающими при этом константами дает требуемое неравенство  \eqref{iq:ufu} с постоянными $C, \overline C_M$.  
\end{proof}

\begin{example}[{на основе одного примера Т.~Лионса [5, XI\,B\,2]}] 
Пусть $\lambda$ --- мера Лебега в единичном шаре $\BB\subset \CC^n$, нормированная условием $\lambda (\BB)=1$. 
Тогда 
\begin{equation}\label{vB}
V_{\lambda}(z):=\int h_{2n}(|z-z'|)\dd \lambda(z')-
h_{2n}(|z|), \quad z\in \BB, 
\end{equation} 
--- {\it положительная тестовая функция для\/ $\BB$ вне\/} 
шара $\varepsilon \BB$ при $\varepsilon \in (0,1)$ [4, предложение 1.4]. Пусть $0\neq z\in \BB$ и $B(z,r):=z+r\BB$ --- шар радиуса $r\in (0,1-|z|)$ с центром $z$,  $\lambda_{z,r}$ --- сужение $\lambda$ на $B(z,r)$. Рассмотрим вероятностную меру $\mu_{z,r}:=\lambda-\lambda_{z,r} +\lambda \bigl(B(z,r)\bigr)\delta_z\geq 0$, где $\delta_z$ --- мера Дирака в точке $z$. Тогда функция $V_{\mu_{z,r}}$, построенная по правилу \eqref{vB} с заменой $\lambda$ на $\mu_{z,r}$, --- {\it вещественная тестовая функция для $\BB$ вне $\varepsilon \BB$\/} при $\varepsilon\in (0,|z|)$ с $V_{\mu_{z,r}}(z)=-\infty$, т.\,е. строго меньшая нуля в открытом множестве, содержащем точку $z$. Таким образом, класс функций $\sbh_0^{\pm}(\BB\setminus \varepsilon \BB; <+\infty)$ строго шире класса $\sbh_0^{+}(\BB\setminus \varepsilon \BB; <+\infty)$ при $\varepsilon \in (0,|z|)$.
\end{example}

\paragraph{Случай $n=1$.} Для области $D\subset \CC_{\infty}:=\CC_{\infty}^1$ теореме \ref{th1} можно придать форму критерия. 
\begin{theorem}\label{th2} Пусть $D$  односвязная и на $\partial D$ более одной точки,  $-\infty \not\equiv M\in \sbh (D)$ --- 
непрерывная функция с мерой Рисса $\nu_M$ и ${\sf Z}:=\{{\sf z}_k\}_{k=1,2,\dots}$ --- последовательность точек в $D$ без предельных точек в $D$.  
Тогда эквивалентны три утверждения:
\begin{enumerate}[{\rm (i)}]
\item\label{a} ${\sf Z}$ --- последовательность нулей для некоторой функции $f\in \Hol(D)$ в том смысле, что для каждой точки $z\in D$ число повторений точки $z$ в последовательности $\sf Z$ равно кратности нуля функции $f$ в точке $z$, 
с ограничением $|f|\leq \exp M$ на $D$;  
\item\label{b} для любых $S\overset{\eqref{DSD}}{\Subset}\Int S_0\subset S_0\subset D$ и $b\in \RR^+\setminus \{0\}$ найдется такое число $C\in \RR^+$, что
\begin{equation}\label{supuM}
 \sum_{{\sf z}_k\in D\setminus S} v({\sf z}_k) \leq \int_{D\setminus S} v \dd \nu_M+C
\quad\text{для всех  $v\in \sbh_0^{\pm}(D\setminus S, S_0;\leq b)$};
\end{equation}

\item\label{c} существуют $S\overset{\eqref{DSD}}{\Subset}\Int S_0\subset S_0\subset D$, $b\in \RR^+\setminus \{0\}$ и $C\in \RR^+$, для которых 
имеем \eqref{supuM}.
\end{enumerate}
\end{theorem}
Импликация \eqref{a}$\Rightarrow$\eqref{b} сразу следует из теоремы \ref{th1}; \eqref{b}$\Rightarrow$\eqref{c} --- очевидно. Доказательство импликации \eqref{c}$\Rightarrow$\eqref{a} использует [6, \S~8, теорема 8.2], [7, 9.3] и будет изложено в ином месте.

\begin{center}
{\sc Литература}
\end{center}

[{\bf 1}] Б.\,Н.~Хабибуллин, А.\,П.~Розит, Функц. анализ и его прил., {\bf 52}:1 (2018), 26--42. \; [{\bf 2}] У.~Хейман, П.~Кеннеди, {\it Субгармонические функции,\/} Мир, М., 1980. \; [{\bf 3}] S.\,J.~Gardiner, {\it Harmonic Approximation,\/} Cambridge Univ. Press, Cambridge, 1995. \; [{\bf 4}] Б.\,Н.~Хабибуллин,
Сиб. матем. журн., {\bf 44}:4 (2003), 905--925.
\; [{\bf 5}] P.~Koosis, {The logarithmic integral. II,\/} 
Cambridge Univ. Press, Cambridge, 1992. \; [{\bf 6}]
Б.\,Н.~Хабибуллин,  Изв. РАН, cер. матем., {\bf 65}:5 (2001), 167--190. 
\; [{\bf 7}]
Б.\,Н.~Хабибуллин, А.\,П.~Розит, Э.\,Б.~Хабибуллина, 
Итоги науки и техн. Сер. Соврем. мат. и ее прил. Темат. обз. (2018), в печати.
\end{document}